\documentclass[11pt]{article}
\usepackage{latexsym,amssymb,amsmath,enumerate,geometry,float,cite,tikz,setspace,verbatim}
\newtheorem{definition}{Definition}
\newtheorem{theorem}[definition]{Theorem}

\newtheorem{lemma}[definition]{Lemma}
\newtheorem{proposition}[definition]{Proposition}

\newtheorem{claim}{Claim}

\usepackage{lineno}

\begin{document}


\onehalfspace

\title{Identifying Codes in the Complementary Prism of Cycles}

\author{
M\'{a}rcia R. Cappelle,
Erika M.M. Coelho,
Hebert Coelho,
Lucia D. Penso,
Dieter Rautenbach}

\date{}

\maketitle

\begin{center}
{\small 
$^1$ Instituto de Inform\'{a}tica, Universidade Federal de Goi\'{a}s, Goi\^{a}nia, Brazil\\
\texttt{marcia@inf.ufg.br, 
erikamorais@inf.ufg.br,
hebert@inf.ufg.br}\\[3mm]
$^2$ Institute of Optimization and Operations Research, Ulm University, Ulm, Germany\\
\texttt{lucia.penso@uni-ulm.de, dieter.rautenbach@uni-ulm.de}
}
\end{center}

\begin{abstract}
We show that an identifying code of minimum order 
in the complementary prism of a cycle of order $n$ 
has order $7n/9+\Theta(1)$.
Furthermore, 
we observe that the clique-width of the complementary prism 
of a graph of clique-width $k$ is at most $4k$,
and discuss some algorithmic consequences.
\end{abstract}

{\small 

\medskip


}

\section{Introduction}\label{section1}

We consider finite, simple, and undirected graphs, and use standard notation and terminology.

For a positive integer $d$, a graph $G$, and a vertex $u$ of $G$, 
let $N^{\leq d}_G[u]$ be the set of vertices of $G$ at distance at most $d$ from $u$.
Note that the closed neighborhood $N_G[u]$ of $u$ in $G$ coincides with $N^{\leq 1}_G[u]$.
A set $C$ of vertices of a graph $G$ is a {\it $d$-identifying code} in $G$ for a positive integer $d$ \cite{kcl}
if the sets $N^{\leq d}_G[u]\cap C$ are non-empty and distinct for all vertices $u$ of $G$.
A $1$-identifying code is known simply as an {\it identifying code}.
Let ${\rm ic}(G)$ denote the minimum order of an identifying code in $C$.

It is algorithmically hard \cite{chl,a} to determine identifying codes of minimum order 
even for planar graphs of arbitrarily large girth.
Exact values, density results, as well as good upper and lower bounds 
have been studied in detail for many special graphs; 
in particular for graphs that arise by product operations using simple factors
such as grids \cite{bchl,bcmms,bl,cghlmpz,chlz,dgm,gms1,gms2,gw,jl,ms}.
In the present paper we study identifying codes in the complementary prism of cycles.
The related notion of locating-domination was studied for such graphs in \cite{hhks}.

Complementary prisms were introduced by Haynes et al. \cite{hhsv} 
as a variation of the well-known {\it prism} of a graph \cite{hik}.
For a graph $G$ with vertex set $V(G)=\{ v_1,\ldots,v_n\}$ and edge set $E(G)$,
the {\it complementary prism of $G$} is the graph denoted by $G\bar{G}$
with vertex set $V(G\bar{G})=\{ v_1,\ldots,v_n\}\cup\{ \bar{v}_1,\ldots,\bar{v}_n\}$
and edge set
\begin{eqnarray*}
E(G\bar{G})&=& E(G)\cup \{ \bar{v}_i\bar{v}_j:1\leq i<j\leq n\mbox{ and }v_iv_j\not\in E(G)\}
\cup \{ v_1\bar{v}_1,\ldots,v_n\bar{v}_n\}.
\end{eqnarray*}
In other words, 
the complementary prism $G\bar{G}$ of $G$ arises from the disjoint union of the graph $G$ and its complement $\bar{G}$
by adding the edges of a perfect matching joining corresponding vertices of $G$ and $\bar{G}$.
For every vertex $u$ of $G$, we will consistently denote the corresponding vertex of $\bar{G}$ by $\bar{u}$,
that is, $V(G\bar{G})=V(G)\cup V(\bar{G})$ where $V(\bar{G})=\{ \bar{v}_1,\ldots,\bar{v}_n\}$.
For a positive integer $k$, let $[k]$ denote the set of positive integers at most $k$.
For an integer $n$ at least $3$, let $C_n$ denote the cycle of order $n$.

In Section \ref{section2} we determine the minimum order of an identifying code in $C_n\bar{C}_n$ up to a small constant.
Note that for $n\geq 6$ and $d\geq 2$, 
the graph $C_n\bar{C}_n$ contains distinct vertices $u$ and $v$ with $N_{C_n\bar{C}_n}^{\leq d}[u]=N_{C_n\bar{C}_n}^{\leq d}[v]$,
which implies that there is no $d$-identifying code in $C_n\bar{C}_n$ for such values.

Before we proceed to Section \ref{section2}, we make some more general algorithmic observations.
In \cite{a} Auger describes an involved linear time dynamic programming algorithm 
that determines an identifying code of minimum order for a given tree. 
In \cite{cghlmm} Charon et al. present a similar algorithm for oriented trees,
and explicitly mention that it is an open issue whether, for any fixed $d$ at least $2$,
it is possible to determine a $d$-identifying code of minimum order for a given tree in polynomial time.
In fact, the existence of such efficient algorithms follows immediately 
from general results \cite{cmr} concerning graph of bounded clique-width, 
such as trees, which have clique-width at most $3$.
For a positive integer $d$, and two vertices $u$ and $v$ of a graph $G$, 
we have $v\in N^{\leq d}_G[u]$ if and only if
\begin{eqnarray*}
\exists v_0,v_1,\ldots,v_d\in V(G) &:& (u=v_0)\wedge (v=v_d)\\
&& \wedge \Big((v_0v_1\in E(G))\vee (v_0=v_1)\Big)
\wedge \cdots
\wedge \Big((v_{d-1}v_d\in E(G))\vee (v_{d-1}=v_d)\Big).
\end{eqnarray*}
Furthermore, a set $C$ of vertices of $G$ is a $d$-identifying code in $G$ if and only if 
\begin{eqnarray*}
\Big(\forall u\in V(G) &:& \exists v\in C:v\in N_G^{\leq d}[u]\Big)\wedge\\
\Big(\forall x,y\in V(G) &:& (x\not=y)\Rightarrow \\
&& \exists z\in C: 
\left((z\in N_G^{\leq d}[x])\wedge (z\not\in N_G^{\leq d}[y])\right)
\vee
\left((z\not\in N_G^{\leq d}[x])\wedge (z\in N_G^{\leq d}[y])\right)\Big).
\end{eqnarray*}
These observations imply that the optimization problem to determine a $d$-identifying code of minimum order
is expressible in the {\it LinEMSOL$(\tau_1)$} logic \cite{cmr}.
Therefore, if $cw$ is some constant, 
and ${\cal G}$ is a class of graphs 
such that every graph $G$ in ${\cal G}$ 
has clique-width at most $cw$, and a clique-width expression  for $G$ 
using at most $cw$ distinct labels 
can be determined in polynomial time,
then $d$-identifying codes of minimum order can be determined in polynomial time for the graphs in ${\cal G}$
(cf. Theorem 4 in \cite{cmr}).
For the class of trees, this immediately implies the existence of linear time algorithms 
that determine a $d$-identifying code of minimum order for any fixed $d$.
These algorithmic consequences extend to complementary prisms by the following result.

\begin{proposition}\label{proposition1}
If $G$ is a graph of clique-width $cw$, then $G\bar{G}$ has clique-width at most $4cw$.
\end{proposition}
{\it Proof:} Let $G$ be a graph of clique-width $cw$.
In \cite{lr} it is shown that there is a rooted binary tree $T$ 
whose leaves are the vertices of $G$
such that, 
for every vertex $s$ of $T$, 
the set $V_s$ of vertices of $G$ that are descendants of $s$ in $T$
partitions into at most $cw$ equivalence classes 
with respect to the equivalence relation $\sim$,
where $u\sim v$ for $u,v\in V_s$ if and only if $N_G[u]\setminus V_s=N_G[v]\setminus V_s$.
Replacing in $T$ every leaf $u$ with parent $x$ by three vertices
$u$, $\bar{u}$, and $y$,
and adding the arcs $(x,y)$, $(y,u)$, and $(y,\bar{u})$,
we obtain a rooted binary tree $T'$ whose leaves are the vertices of $G\bar{G}$.
By the definition of $G\bar{G}$, we obtain that 
for every vertex $s'$ of $T'$, 
the set $V'_s$ of vertices of $G\bar{G}$ that are descendants of $s'$ in $T'$
partitions into at most $2cw$ equivalence classes 
with respect to the equivalence relation $\sim'$,
where $u\sim' v$ for $u,v\in V'_s$ if and only if $N_{G\bar{G}}[u]\setminus V'_s=N_{G\bar{G}}[v]\setminus V'_s$.
Again by \cite{lr}, this implies that the clique-width of $G\bar{G}$ is at most $4cw$. $\Box$

\section{Minimum identifying code in $C_n\bar{C}_n$}\label{section2}

Throughout this section, let $C_n:v_1v_2\ldots v_nv_1$ be a cycle of order $n$ at least $3$,
and let $G=C_n\bar{C}_n$. We identify indices of vertices of $G$ modulo $n$.
For a subset $C$ of $V(C_n)$, let $x(C)$ denote the characteristic vector of $C$, 
that is, $x(C)=(x_1,\ldots,x_n)\in \{ 0,1\}^n$ where $x_i=1$ if and only if $v_i\in C$ for $i\in [n]$.
Similarly, for a subset $\bar{C}$ of $V(\bar{C}_n)$, 
let $x(\bar{C})=(\bar{x}_1,\ldots,\bar{x}_n)\in \{ 0,1\}^n$ where $\bar{x}_i=1$ 
if and only if $\bar{v}_i\in \bar{C}$ for $i\in [n]$.

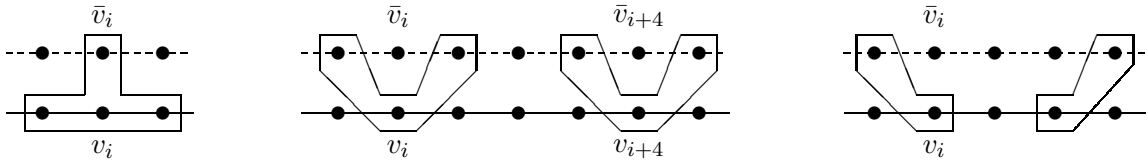
\begin{figure}[H]
\begin{center}
$\mbox{}$ \hfill 
\unitlength 0.8mm 
\linethickness{0.4pt}
\ifx\plotpoint\undefined\newsavebox{\plotpoint}\fi 
\begin{picture}(35,26)(0,0)
\put(10,10){\circle*{2}}
\put(20,10){\circle*{2}}
\put(30,10){\circle*{2}}
\put(10,20){\circle*{2}}
\put(20,20){\circle*{2}}
\put(30,20){\circle*{2}}
\put(17,23){\line(1,0){6}}
\put(23,23){\line(0,-1){10}}
\put(23,13){\line(1,0){10}}
\put(33,13){\line(0,-1){6}}
\put(33,7){\line(-1,0){26}}
\put(7,7){\line(0,1){6}}
\put(7,13){\line(1,0){10}}
\put(17,13){\line(0,1){10}}
\put(20,4){\makebox(0,0)[cc]{$v_i$}}
\put(20,26){\makebox(0,0)[cc]{$\bar{v}_i$}}
\put(4,10){\line(1,0){31}}
\put(3.93,19.93){\line(1,0){.9688}}
\put(5.867,19.93){\line(1,0){.9688}}
\put(7.805,19.93){\line(1,0){.9688}}
\put(9.742,19.93){\line(1,0){.9688}}
\put(11.68,19.93){\line(1,0){.9688}}
\put(13.617,19.93){\line(1,0){.9688}}
\put(15.555,19.93){\line(1,0){.9687}}
\put(17.492,19.93){\line(1,0){.9688}}
\put(19.43,19.93){\line(1,0){.9688}}
\put(21.367,19.93){\line(1,0){.9688}}
\put(23.305,19.93){\line(1,0){.9688}}
\put(25.242,19.93){\line(1,0){.9688}}
\put(27.18,19.93){\line(1,0){.9688}}
\put(29.117,19.93){\line(1,0){.9688}}
\put(31.055,19.93){\line(1,0){.9688}}
\put(32.992,19.93){\line(1,0){.9688}}
\end{picture}
\hfill
\linethickness{0.4pt}
\ifx\plotpoint\undefined\newsavebox{\plotpoint}\fi 
\begin{picture}(75,26)(0,0)
\put(10,10){\circle*{2}}
\put(50,10){\circle*{2}}
\put(20,10){\circle*{2}}
\put(60,10){\circle*{2}}
\put(30,10){\circle*{2}}
\put(70,10){\circle*{2}}
\put(40,10){\circle*{2}}
\put(10,20){\circle*{2}}
\put(50,20){\circle*{2}}
\put(20,20){\circle*{2}}
\put(60,20){\circle*{2}}
\put(30,20){\circle*{2}}
\put(70,20){\circle*{2}}
\put(40,20){\circle*{2}}
\put(20,4){\makebox(0,0)[cc]{$v_i$}}
\put(60,4){\makebox(0,0)[cc]{$v_{i+4}$}}
\put(20,26){\makebox(0,0)[cc]{$\bar{v}_i$}}
\put(60,26){\makebox(0,0)[cc]{$\bar{v}_{i+4}$}}
\put(7,23){\line(1,0){6}}
\put(47,23){\line(1,0){6}}
\put(13,23){\line(2,-5){4}}
\put(53,23){\line(2,-5){4}}
\put(17,13){\line(1,0){6}}
\put(57,13){\line(1,0){6}}
\put(23,13){\line(2,5){4}}
\put(63,13){\line(2,5){4}}
\put(27,23){\line(1,0){6}}
\put(67,23){\line(1,0){6}}
\put(33,23){\line(0,-1){6}}
\put(73,23){\line(0,-1){6}}
\put(33,17){\line(-1,-1){10}}
\put(73,17){\line(-1,-1){10}}
\put(23,7){\line(-1,0){6}}
\put(63,7){\line(-1,0){6}}
\put(17,7){\line(-1,1){10}}
\put(57,7){\line(-1,1){10}}
\put(7,17){\line(0,1){6}}
\put(47,17){\line(0,1){6}}
\put(3.93,19.93){\line(1,0){.9861}}
\put(5.902,19.93){\line(1,0){.9861}}
\put(7.874,19.93){\line(1,0){.9861}}
\put(9.846,19.93){\line(1,0){.9861}}
\put(11.819,19.93){\line(1,0){.9861}}
\put(13.791,19.93){\line(1,0){.9861}}
\put(15.763,19.93){\line(1,0){.9861}}
\put(17.735,19.93){\line(1,0){.9861}}
\put(19.707,19.93){\line(1,0){.9861}}
\put(21.68,19.93){\line(1,0){.9861}}
\put(23.652,19.93){\line(1,0){.9861}}
\put(25.624,19.93){\line(1,0){.9861}}
\put(27.596,19.93){\line(1,0){.9861}}
\put(29.569,19.93){\line(1,0){.9861}}
\put(31.541,19.93){\line(1,0){.9861}}
\put(33.513,19.93){\line(1,0){.9861}}
\put(35.485,19.93){\line(1,0){.9861}}
\put(37.457,19.93){\line(1,0){.9861}}
\put(39.43,19.93){\line(1,0){.9861}}
\put(41.402,19.93){\line(1,0){.9861}}
\put(43.374,19.93){\line(1,0){.9861}}
\put(45.346,19.93){\line(1,0){.9861}}
\put(47.319,19.93){\line(1,0){.9861}}
\put(49.291,19.93){\line(1,0){.9861}}
\put(51.263,19.93){\line(1,0){.9861}}
\put(53.235,19.93){\line(1,0){.9861}}
\put(55.207,19.93){\line(1,0){.9861}}
\put(57.18,19.93){\line(1,0){.9861}}
\put(59.152,19.93){\line(1,0){.9861}}
\put(61.124,19.93){\line(1,0){.9861}}
\put(63.096,19.93){\line(1,0){.9861}}
\put(65.069,19.93){\line(1,0){.9861}}
\put(67.041,19.93){\line(1,0){.9861}}
\put(69.013,19.93){\line(1,0){.9861}}
\put(70.985,19.93){\line(1,0){.9861}}
\put(72.957,19.93){\line(1,0){.9861}}
\put(4,10){\line(1,0){71}}
\end{picture}
\hfill
\linethickness{0.4pt}
\ifx\plotpoint\undefined\newsavebox{\plotpoint}\fi 
\begin{picture}(55,26)(0,0)
\put(10,10){\circle*{2}}
\put(30,10){\circle*{2}}
\put(50,10){\circle*{2}}
\put(20,10){\circle*{2}}
\put(40,10){\circle*{2}}
\put(10,20){\circle*{2}}
\put(30,20){\circle*{2}}
\put(50,20){\circle*{2}}
\put(20,20){\circle*{2}}
\put(40,20){\circle*{2}}
\put(20,4){\makebox(0,0)[cc]{$v_i$}}
\put(20,26){\makebox(0,0)[cc]{$\bar{v}_i$}}
\put(7,23){\line(1,0){6}}
\put(13,23){\line(2,-5){4}}
\put(17,7){\line(-1,1){10}}
\put(7,17){\line(0,1){6}}
\put(17,13){\line(1,0){6}}
\put(23,13){\line(0,-1){6}}
\put(37,7){\line(0,1){6}}
\put(37,13){\line(1,0){6}}
\put(43,13){\line(2,5){4}}
\put(47,23){\line(1,0){6}}
\put(53,23){\line(0,-1){5}}
\put(37,7){\line(1,0){6}}
\multiput(43,7)(.0336700337,.037037037){297}{\line(0,1){.037037037}}
\put(17,7){\line(1,0){6}}
\put(5,10){\line(1,0){50}}
\put(4.93,19.93){\line(1,0){.9804}}
\put(6.89,19.93){\line(1,0){.9804}}
\put(8.851,19.93){\line(1,0){.9804}}
\put(10.812,19.93){\line(1,0){.9804}}
\put(12.773,19.93){\line(1,0){.9804}}
\put(14.734,19.93){\line(1,0){.9804}}
\put(16.694,19.93){\line(1,0){.9804}}
\put(18.655,19.93){\line(1,0){.9804}}
\put(20.616,19.93){\line(1,0){.9804}}
\put(22.577,19.93){\line(1,0){.9804}}
\put(24.538,19.93){\line(1,0){.9804}}
\put(26.498,19.93){\line(1,0){.9804}}
\put(28.459,19.93){\line(1,0){.9804}}
\put(30.42,19.93){\line(1,0){.9804}}
\put(32.381,19.93){\line(1,0){.9804}}
\put(34.341,19.93){\line(1,0){.9804}}
\put(36.302,19.93){\line(1,0){.9804}}
\put(38.263,19.93){\line(1,0){.9804}}
\put(40.224,19.93){\line(1,0){.9804}}
\put(42.185,19.93){\line(1,0){.9804}}
\put(44.145,19.93){\line(1,0){.9804}}
\put(46.106,19.93){\line(1,0){.9804}}
\put(48.067,19.93){\line(1,0){.9804}}
\put(50.028,19.93){\line(1,0){.9804}}
\put(51.989,19.93){\line(1,0){.9804}}
\put(53.949,19.93){\line(1,0){.9804}}
\end{picture}
\hfill $\mbox{}$
\end{center}
\caption{Condition $C(i)$ implies that at least one of the four vertices indicated in the left figure belongs to $C\cup \bar{C}$.
Condition $\bar{C}(i,i+4)$ implies that at least one of the six vertices indicated in the middle figure belongs to $C\cup \bar{C}$.
Condition $\bar{C}(i,i+2)$ implies that at least one of the four vertices indicated in the right figure belongs to $C\cup \bar{C}$.
Note that for $\bar{C}_n$, instead of indicating the edges, we indicate the non-edges by dashed lines.}\label{fig1}
\end{figure}

\begin{lemma}\label{lemma1}
For an integer $n$ at least $9$, let $G=C_n\bar{C}_n$.
Let $C\subseteq V(C_n)$ and $\bar{C}\subseteq V(\bar{C}_n)$.
Let $x(C)=(x_1,\ldots,x_n)$ and $x(\bar{C})=(\bar{x}_1,\ldots,\bar{x}_n)$.

If $C\cup \bar{C}$ is an identifying code in $G$, 
then the following conditions hold for every $i,j\in [n]$ with $(j-i)\,\,{\rm mod}\,\, n\not\in \{ 0,2\}$ (cf. Figure \ref{fig1}):
$$
\begin{array}{lcll}
C(i) &:& x_{i-1}+x_i+\bar{x}_i+x_{i+1} & \geq 1,\\
C(i,i+1) &:& x_{i-1}+\bar{x}_i+\bar{x}_{i+1}+x_{i+2} & \geq 1,\\
C(i,i+2) &:& x_{i-1}+x_i+\bar{x}_i+x_{i+2}+\bar{x}_{i+2}+x_{i+3} & \geq 1,\\
\bar{C}(i,j) &:& \bar{x}_{i-1}+x_i+\bar{x}_{i+1}+\bar{x}_{j-1}+x_j+\bar{x}_{j+1} & \geq 1,\mbox{ and}\\
\bar{C}(i,i+2) &:& \bar{x}_{i-1}+x_i+x_{i+2}+\bar{x}_{i+3} & \geq 1.
\end{array}
$$
Furthermore, if $|\bar{C}|\geq 4$, then $C\cup \bar{C}$ is an identifying code in $G$
if and only if these conditions hold.
\end{lemma}
{\it Proof:} Note that for distinct vertices $u$ and $v$ of $G$, 
we have $N_G[u]\cap (C\cup \bar{C})\not=N_G[v]\cap (C\cup \bar{C})$
if and only if $C\cup \bar{C}$ intersects $(N_G[u]\setminus N_G[v])\cup (N_G[v]\setminus N_G[u])$.
Therefore, for $i\in [n]$, we have that
$C(i)$ is equivalent to $N_G(v_i)\cap (C\cup \bar{C})\not=\emptyset$,
$C(i,i+1)$ is equivalent to $N_G(v_i)\cap (C\cup \bar{C})\not=N_G(v_{i+1})\cap (C\cup \bar{C})$,
$C(i,i+2)$ is equivalent to $N_G(v_i)\cap (C\cup \bar{C})\not=N_G(v_{i+2})\cap (C\cup \bar{C})$, and
$\bar{C}(i,i+2)$ is equivalent to $N_G(\bar{v}_i)\cap (C\cup \bar{C})\not=N_G(\bar{v}_{i+2})\cap (C\cup \bar{C})$.
For $i,j\in [n]$ with $(j-i)\,\,{\rm mod}\,\, n\geq 3$, we have that
$C(i)$ and $C(j)$ together are equivalent to $N_G(v_i)\cap (C\cup \bar{C})\not=N_G(v_j)\cap (C\cup \bar{C})$.
For $i,j\in [n]$ with $(j-i)\,\,{\rm mod}\,\, n\not\in \{ 0,2\}$, we have that
$\bar{C}(i,j)$ is equivalent to $N_G(\bar{v}_i)\cap (C\cup \bar{C})\not=N_G(\bar{v}_j)\cap (C\cup \bar{C})$.
Hence, all these conditions are necessary.
Note that $|\bar{C}|\geq 4$ implies 
$N_G(v_i)\cap (C\cup \bar{C})\not=N_G(\bar{v}_j)\cap (C\cup \bar{C})\not=\emptyset$ for every $i,j\in [n]$,
in which case, the given conditions are also sufficient. $\Box$

\begin{figure}[H]
\begin{center}
\unitlength 0.7mm 
\linethickness{0.4pt}
\ifx\plotpoint\undefined\newsavebox{\plotpoint}\fi 
\begin{picture}(230,18)(0,0)
\put(4,4){\circle*{2}}
\put(112,4){\circle*{2}}
\put(202,4){\circle*{2}}
\put(34,4){\circle*{2}}
\put(142,4){\circle*{2}}
\put(64,4){\circle*{2}}
\put(172,4){\circle*{2}}
\put(14,4){\circle*{2}}
\put(122,4){\circle*{2}}
\put(212,4){\circle*{2}}
\put(44,4){\circle*{2}}
\put(152,4){\circle*{2}}
\put(74,4){\circle*{2}}
\put(182,4){\circle*{2}}
\put(24,4){\circle*{2}}
\put(132,4){\circle*{2}}
\put(222,4){\circle*{2}}
\put(54,4){\circle*{2}}
\put(162,4){\circle*{2}}
\put(84,4){\circle*{2}}
\put(192,4){\circle*{2}}
\put(4,14){\circle*{2}}
\put(112,14){\circle*{2}}
\put(202,14){\circle*{2}}
\put(34,14){\circle*{2}}
\put(142,14){\circle*{2}}
\put(64,14){\circle*{2}}
\put(172,14){\circle*{2}}
\put(14,14){\circle*{2}}
\put(122,14){\circle*{2}}
\put(212,14){\circle*{2}}
\put(44,14){\circle*{2}}
\put(152,14){\circle*{2}}
\put(74,14){\circle*{2}}
\put(182,14){\circle*{2}}
\put(24,14){\circle*{2}}
\put(132,14){\circle*{2}}
\put(222,14){\circle*{2}}
\put(54,14){\circle*{2}}
\put(162,14){\circle*{2}}
\put(84,14){\circle*{2}}
\put(192,14){\circle*{2}}
\put(1,1){\framebox(26,6)[cc]{}}
\put(109,1){\framebox(26,6)[cc]{}}
\put(199,1){\framebox(26,6)[cc]{}}
\put(40,11){\framebox(37,6)[cc]{}}
\put(148,11){\framebox(37,6)[cc]{}}
\put(0,0){\framebox(88,18)[cc]{}}
\put(108,0){\framebox(88,18)[cc]{}}
\put(230,10){\makebox(0,0)[cc]{$\cdots$}}
\put(98,9){\makebox(0,0)[cc]{$\cdots$}}
\end{picture}
\end{center}
\caption{Some identifying code for $C_n\bar{C_n}$.}\label{fig2}
\end{figure}
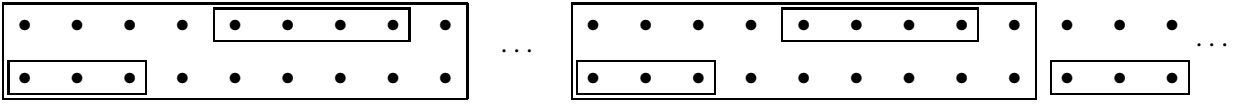

\begin{lemma}\label{lemma2}
For an integer $n$ at least $9$, let $G=C_n\bar{C}_n$.
Let $C\subseteq V(C_n)$ and $\bar{C}\subseteq V(\bar{C}_n)$.
Let $x(C)=(x_1,\ldots,x_n)$ and $x(\bar{C})=(\bar{x}_1,\ldots,\bar{x}_n)$.

If $k=\left\lfloor\frac{n}{9}\right\rfloor$, and for $i\in [n]$ (cf. Figure \ref{fig2}),
$$
x_i=
\left\{
\begin{array}{ll}
1, & i\,\,{\rm mod}\,\, 9\in \{ 1,2,3\}\mbox{ and }i\leq 9k,\\
1, & i\geq 9k+1,\mbox{ and }\\ 
0, & \mbox{otherwise}
\end{array}
\right.
$$
and
$$
\bar{x}_i=
\left\{
\begin{array}{ll}
1, & i\,\,{\rm mod}\,\, 9\in \{ 5,6,7,8\}\mbox{ and }i\leq 9k,\mbox{ and }\\ 
0, & \mbox{otherwise},
\end{array}
\right.
$$
then $C\cup \bar{C}$ is an identifying code in $G$.
In particular,
$${\rm ic}(G)\leq \frac{7}{9}n+\frac{16}{9}.$$
\end{lemma}
{\it Proof:} Lemma \ref{lemma1} easily implies that $C\cup \bar{C}$ is an identifying code in $G$.
Furthermore,
$|C\cup \bar{C}|
=7k+(n-9k)
=n-2k
=n-2\left\lfloor\frac{n}{9}\right\rfloor
\leq n-2\frac{(n-8)}{9}
=\frac{7}{9}n+\frac{16}{9}$.
$\Box$

\begin{lemma}\label{lemma2b}
For an integer $n$ at least $9$, let $G=C_n\bar{C}_n$.
Let $C\subseteq V(C_n)$ and $\bar{C}\subseteq V(\bar{C}_n)$
be such that $C\cup \bar{C}$ is an identifying code in $G$.
Let $x(C)=(x_1,\ldots,x_n)$ and $x(\bar{C})=(\bar{x}_1,\ldots,\bar{x}_n)$.
Let $\bar{I}=\{ i\in [n]:\bar{x}_{i-1}+x_i+\bar{x}_{i+1}=0\}$.

If $|\bar{C}|\geq 6$ and $i\in [n]$ is such that $i-5,i,i+1,i+6\not\in \bar{I}$ and 
${\bar{x}_i\,\,\bar{x}_{i+1}\choose x_i\,\,x_{i+1}}={0\,\,0\choose 0\,\,0}$, then
\begin{enumerate}[(i)]
\item either there are subsets $C'\subseteq V(C_n)$ and $\bar{C}'\subseteq V(\bar{C}_n)$ such that 
\begin{itemize}
\item $C'\cup \bar{C}'$ is an identifying code in $G$ with $|C'\cup \bar{C}'|\leq |C\cup \bar{C}|$, and
\item 
$$\left|\left\{ j\in[n]:{\bar{x}'_j\choose x'_j}={0\choose 0}\right\}\right|
<\left|\left\{ j\in[n]:{\bar{x}_j\choose x_j}={0\choose 0}\right\}\right|
$$
where $x(C')=(x'_1,\ldots,x'_n)$ and $x(\bar{C}')=(\bar{x}'_1,\ldots,\bar{x}'_n)$.
\end{itemize}
\item or 
$\left\{ {\bar{x}_{i-1}\,\,\bar{x}_i\,\,\ldots\,\,\bar{x}_{i+7}\choose x_{i-1}\,\,x_i\,\,\ldots\,\,x_{i+7}},
{\bar{x}_{i-6}\,\,\bar{x}_{i-5}\,\,\ldots\,\,\bar{x}_{i+2}\choose x_{i-6}\,\,x_{i-5}\,\,\ldots\,\,x_{i+2}}\right\}$
contains 
${1\,\,0\,\,0\,\,1\,\,0\,\,1\,\,0\,\,0\,\,1\choose 1\,\,0\,\,0\,\,1\,\,0\,\,1\,\,0\,\,0\,\,1}$.
\end{enumerate}
\end{lemma}
{\it Proof:} We use the conditions from Lemma \ref{lemma1}.
Let $i\in [n]$ be such that $i-5,i,i+1,i+6\not\in \bar{I}$ and 
${\bar{x}_i\,\,\bar{x}_{i+1}\choose x_i\,\,x_{i+1}}={0\,\,0\choose 0\,\,0}$.
$C(i)$ and $C(i+1)$ imply $x_{i-1}=x_{i+2}=1$.
Since $i,i+1\not\in \bar{I}$, we have $\bar{x}_{i-1}=\bar{x}_{i+2}=1$.

Let $C'=C\cup \{ v_i,v_{i+1}\}$ and $\bar{C}'=\bar{C}\setminus \{ \bar{v}_{i-1},\bar{v}_{i+2}\}$.
Let $x(C')=(x'_1,\ldots,x'_n)$ and $x(\bar{C}')=(\bar{x}'_1,\ldots,\bar{x}'_n)$.
If $C'\cup \bar{C}'$ is an identifying code in $G$, then (i) holds.
Hence, we may assume that $C'\cup \bar{C}'$ is not an identifying code in $G$.
Since $|\bar{C}'|\geq 4$, some condition from Lemma \ref{lemma1} is violated by $C'\cup \bar{C}'$.
Since $(C\cup \bar{C})\setminus (C'\cup \bar{C}')=\{ \bar{v}_{i-1},\bar{v}_{i+2}\}$, 
a violated condition must involve $\bar{x}'_{i-1}$ or $\bar{x}'_{i+2}$.
By symmetry, we may assume that $\bar{x}'_{i+2}$ is involved in a violated condition.
The conditions that involve $\bar{x}'_{i+2}$ are
$C(i+2)$,
$C(i+1,i+2)$,
$C(i+2,i+3)$,
$C(i,i+2)$,
$C(i+2,i+4)$,
$\bar{C}(i+1,j)$,
$\bar{C}(i+3,j)$,
$\bar{C}(i-1,i+1)$, and
$\bar{C}(i+3,i+5)$ 
for $j\in [n]$ with $(j-i)\,\,{\rm mod}\,\, n\not\in \{ 0,2\}$,
where we replace $x_j$ with $x_j'$ and $\bar{x}_j$ with $\bar{x}'_j$ for all $j\in [n]$.
Since $x'_i=1$,
the conditions 
$C(i+1,i+2)$ and
$C(i,i+2)$
are not violated.
Since $x'_{i+1}=1$,
the conditions 
$C(i+2,i+3)$,
$\bar{C}(i+1,j)$, and
$\bar{C}(i-1,i+1)$
are not violated.
Since $x'_{i+2}=1$,
the conditions 
$C(i+2)$ and
$C(i+2,i+4)$
are not violated.
If $\bar{C}(i+3,j)$ is violated, 
then $x_{i+3}=x'_{i+3}=0$ and $\bar{x}_{i+4}=\bar{x}'_{i+4}=0$.
If $\bar{C}(i+3,i+5)$ is violated, 
then $x_{i+3}=x_{i+5}=0$ and $\bar{x}_{i+6}=0$.
Therefore, by symmetry, we may assume that 
\begin{itemize}
\item either $x_{i+3}=\bar{x}_{i+4}=0$,
\item or $x_{i+3}=x_{i+5}=\bar{x}_{i+6}=0$, and $\bar{x}_{i+4}=1$.
\end{itemize}
In the first case, $\bar{C}(i+1,i+3)$ is violated.
Hence, we may assume $x_{i+3}=x_{i+5}=\bar{x}_{i+6}=0$, and $\bar{x}_{i+4}=1$.
$C(i+1,i+3)$ implies that $\bar{x}_{i+3}=1$ or $x_{i+4}=1$.
Let $C''=C\setminus \{ v_{i+2}\}$ and $\bar{C}''=\bar{C}\cup\{ \bar{v}_{i+1}\}$.
If $C''\cup \bar{C}''$ is an identifying code in $G$, then (i) holds.
Let $x(C'')=(x''_1,\ldots,x''_n)$.
Hence, we may assume that $C''\cup \bar{C}''$ is not an identifying code in $G$.
Since $|\bar{C}''|\geq 4$, some condition from Lemma \ref{lemma1} is violated by $C''\cup \bar{C}''$.
Since $(C\cup \bar{C})\setminus (C''\cup \bar{C}'')=\{ v_{i+2}\}$, 
a violated condition must involve $x''_{i+2}$.
Arguing as above, Lemma \ref{lemma1} implies that 
$\bar{x}_{i+3}=\bar{x}_{i+5}=x_{i+6}=0$.
As noted above, $\bar{x}_{i+3}=0$ implies $x_{i+4}=1$.
Now $C(i+6)$ and $i+6\not\in \bar{I}$ imply $x_{i+7}=\bar{x}_{i+7}=1$,
that is, (ii) holds, which completes the proof. 
$\Box$

\begin{lemma}\label{lemma4}
If $n$ is an integer at least $9$, then ${\rm ic}(C_n\bar{C_n})\geq \frac{7}{9}n-12$.
\end{lemma}
{\it Proof:} We prove the statement by induction on $n$, and use the conditions from Lemma \ref{lemma1}.
Clearly, we may assume that $n>\left\lfloor\frac{9\cdot 12}{7}\right\rfloor=15$.
Let $G=C_n\bar{C}_n$, and 
let $C\subseteq V(C_n)$ and $\bar{C}\subseteq V(\bar{C}_n)$ be such that 
\begin{itemize}
\item $C\cup \bar{C}$ is an identifying code in $G$ with ${\rm ic}(G)=|C\cup \bar{C}|$, and
\item subject to the previous condition,
$$\left|\left\{ i\in[n]:{\bar{x}_i\choose x_i}={0\choose 0}\right\}\right|$$
is as small as possible. 
\end{itemize}
Let $\bar{I}=\{ i\in [n]:\bar{x}_{i-1}+x_i+\bar{x}_{i+1}=0\}$.
By $\bar{C}(i,j)$ and $\bar{C}(i,i+2)$, 
we may assume that $\bar{I}\subseteq [1]$.

If $|\bar{C}|\leq 5$, then there are at least $n-1-2|\bar{C}|$ indices $i$ with $i\not\in \bar{I}$
and $\bar{x}_{i-1}=\bar{x}_{i+1}=0$, which implies $x_i=1$.
Therefore, $|C|+|\bar{C}|\geq (n-1-2|\bar{C}|)+|\bar{C}|
=n-1-|\bar{C}|\geq n-6>\frac{7}{9}n-12$.
Hence, we may assume that $|\bar{C}|\geq 6$.

\begin{claim}\label{claim2}
There is no integer $i$ with $7\leq i\leq n-6$ such that 
${\bar{x}_i\,\,\bar{x}_{i+1}\choose x_i\,\,x_{i+1}}={0\,\,0\choose 0\,\,0}$.
\end{claim}
{\it Proof of Claim \ref{claim2}:}
If there is some integer $i$ with $7\leq i\leq n-6$ such that 
${\bar{x}_i\,\,\bar{x}_{i+1}\choose x_i\,\,x_{i+1}}={0\,\,0\choose 0\,\,0}$,
then $i-5,i,i+1,i+6\not\in \bar{I}$.
Now Lemma \ref{lemma2b} and the choice of $C\cup \bar{C}$
imply 
${\bar{x}_{i-1}\,\,\bar{x}_i\,\,\ldots\,\,\bar{x}_{i+7}\choose x_{i-1}\,\,x_i\,\,\ldots\,\,x_{i+7}}
={1\,\,0\,\,0\,\,1\,\,0\,\,1\,\,0\,\,0\,\,1\choose 1\,\,0\,\,0\,\,1\,\,0\,\,1\,\,0\,\,0\,\,1}$
or
${\bar{x}_{i-6}\,\,\bar{x}_{i-5}\,\,\ldots\,\,\bar{x}_{i+2}\choose x_{i-6}\,\,x_{i-5}\,\,\ldots\,\,x_{i+2}}
={1\,\,0\,\,0\,\,1\,\,0\,\,1\,\,0\,\,0\,\,1\choose 1\,\,0\,\,0\,\,1\,\,0\,\,1\,\,0\,\,0\,\,1}$.
By symmetry, we may assume that the former case occurs.
Let $C'\subseteq V(C_{n-5})$ and $\bar{C}'\subseteq V(\bar{C}_{n-5})$ be such that 
$$
{
\bar{x}'_1\,\,\ldots\,\,\bar{x}'_{n-5}
\choose 
x'_1\,\,\ldots\,\,x'_{n-5}
}
=
{
\bar{x}_1\,\,\ldots\,\,\bar{x}_{i+2}\,\,\bar{x}_{i+8}\,\,\ldots\,\,\bar{x}_n
\choose
x_1\,\,\ldots\,\,x_{i+2}\,\,x_{i+8}\,\,\ldots\,\,x_n
}
$$
for $x(C')=(x'_1,\ldots,x'_{n-5})$ and $x(\bar{C}')=(\bar{x}'_1,\ldots,\bar{x}'_{n-5})$.
Since $|\bar{C}|\geq 6$, we have $|\bar{C}'|\geq 4$.
Considering the conditions from Lemma \ref{lemma1} easily implies that 
$C'\cup \bar{C}'$ is an identifying code in $C_{n-5}\bar{C}_{n-5}$.
By induction, we obtain
$|C\cup \bar{C}|=|C'\cup \bar{C}'|+4
\geq \frac{7}{9}(n-5)-12+4
> \frac{7}{9}n-12$.
$\Box$

\bigskip

\noindent Let $i_1<i_2<\ldots<i_k$ be the increasing sequence of integers $i$ 
with $8\leq i\leq n-6$ such that ${\bar{x}_i\choose x_i}={0\choose 0}$.
Note that, for $j\in [k]$, we have
${\bar{x}_{i_j-1}\choose x_{i_j-1}},{\bar{x}_{i_j+1}\choose x_{i_j+1}}\not={0\choose 0}$.

For $j\in [k-1]$, let $I_j=\{ i\in [n]:i_j\leq i\leq i_{j+1}-1\}$.
Note that $|I_j|\geq 2$ for $j\in [k-1]$.

For $I\subseteq [n]$, let $V(I)=\bigcup_{i\in I}\{ v_i,\bar{v}_i\}$.

\begin{claim}\label{claim3}
If $k\geq 2$, 
then there are integers $\ell,j_1,j_2,\ldots,j_\ell$ with $\ell\geq 2$ and $1=j_1<j_2<\ldots<j_\ell=k$ such that
\begin{eqnarray}
\left|
\left(C\cup \bar{C}\right)
\cap 
V(I_{j_1}\cup \cdots \cup I_{j_2-1})
\right| &\geq &\frac{7}{9}\left|I_{j_1}\cup \cdots \cup I_{j_2-1}\right|-\frac{5}{9},\label{e1}\\
\left|
\left(C\cup \bar{C}\right)
\cap 
V(I_{j_r}\cup \cdots \cup I_{j_{r+1}-1})
\right|&\geq &\frac{7}{9}\left|I_{j_r}\cup \cdots \cup I_{j_{r+1}-1}\right|,
\mbox{ for $r\in [\ell-2]\setminus [1]$, and}\label{e2}\\
\left|
\left(C\cup \bar{C}\right)
\cap 
V(I_{j_{\ell-1}}\cup \cdots \cup I_{j_\ell-1})
\right|&\geq &\frac{7}{9}\left|I_{j_{\ell-1}}\cup \cdots \cup I_{j_\ell-1}\right|-\frac{5}{9}.\label{e3}
\end{eqnarray}
\end{claim}
{\it Proof of Claim \ref{claim3}:}
If for some $j\in [k-1]$, there is some $i\in I_j$ with ${\bar{x}_i\choose x_i}={1\choose 1}$,
then $I_j$ is {\it dirty}; otherwise $I_j$ is {\it clean}. 
Note that, 
if $I_j$ is dirty, then
\begin{eqnarray}
\left|\left(C\cup \bar{C}\right)\cap V(I_j)\right| & \geq & |I_j|,\label{e4}
\end{eqnarray}
and, 
if $I_j$ is clean, then, since $|I_j|\geq 2$,
\begin{eqnarray}
\left|\left(C\cup \bar{C}\right)\cap V(I_j)\right| &=&|I_j|-1
\geq \frac{7}{9}|I_j|-\frac{5}{9}.\label{e5}
\end{eqnarray}

Let $j_1=1$.

Clearly, ${\bar{x}_{i_1-1}\choose x_{i_1-1}}\not={0\choose 0}$.
If ${\bar{x}_{i_1-1}\choose x_{i_1-1}}\not={1\choose 1}$,
then the definition of $j_2$ follows the pattern of the definition of $j_{r+1}$ for $r\geq 2$ described below,
that is, in this case, (\ref{e2}) will be satisfied also for $r=1$,
which is a stronger inequality.
If ${\bar{x}_{i_1-1}\choose x_{i_1-1}}={1\choose 1}$,
then let $j_2$ be maximum such that $j_1<j_2\leq k$ and
$I_j$ is dirty for $j\in \{ j_1,j_1+1,\ldots,j_2-2\}$.
Note that, if $I_{j_1}$ is clean, then $j_2=j_1+1$.
By (\ref{e4}) and (\ref{e5}), we obtain that (\ref{e1}) holds.
If $j_2=k$, then set $\ell=2$, 
and terminate the definition of the sequence $j_1,\ldots,j_\ell$.
Note that (\ref{e3}) coincides with (\ref{e1}) in this case.
If $j_2<k$, then, by the choice of $j_2$, we have that $I_{j_2-1}$ is clean,
which implies that 
${\bar{x}_{i_{j_2}-1}\choose x_{i_{j_2}-1}}\not={1\choose 1}$.

Therefore, we may now assume that for some non-negative integer $r$,
the indices $1=j_1<\cdots <j_r<k$ have already been defined 
in such a way that the corresponding conditions are satisfied,
and that  
${\bar{x}_{i_{j_r}-1}\choose x_{i_{j_r}-1}}\not\in \left\{ {0\choose 0},{1\choose 1}\right\}$.
We will define $j_{r+1}$ with $j_r<j_{r+1}\leq k$ such that the corresponding condition is satisfied.
We consider different cases.
In each case, we consider potential choices $j'_{r+1}$ and possibly $j''_{r+1}$ for $j_{r+1}$.
As before, if one of $j_{r+1}$, $j'_{r+1}$, or $j''_{r+1}$ equals $k$, 
then set $\ell=r+1$, and terminate the definition of the sequence $j_1,\ldots,j_\ell$.
In such a case, (\ref{e4}) and (\ref{e5}) will imply (\ref{e3}).

Let $t=i_{j_r}$.

\bigskip

\noindent {\bf Case 1} {\it $I_{j_r}$ is clean and ${\bar{x}_t\,\,\bar{x}_{t+1}\choose x_t\,\,x_{t+1}}={0\,\,1\choose 0\,\,0}$.}

\bigskip 

\noindent $C(t)$ implies $x_{t-1}=1$, and hence, $\bar{x}_{t-1}=0$.
Since $t+1\not\in \bar{I}$, we have $\bar{x}_{t+2}=1$.
Since $I_{j_r}$ is clean, $x_{t+2}=0$.
$\bar{C}(t,t+2)$ implies $\bar{x}_{t+3}=1$.
Since $I_{j_r}$ is clean, $x_{t+3}=0$.
$\bar{C}(t+1,t+3)$ implies $\bar{x}_{t+4}=1$.
Since $I_{j_r}$ is clean, $x_{t+4}=0$.
This implies $i_{j_r+1}-i_{j_r}\geq 5$.
Since 
$$\left|\left(C\cup \bar{C}\right)\cap V(I_{j_r})\right|
=|I_{j_r}|-1\geq \frac{4}{5}|I_{j_r}|>\frac{7}{9}|I_{j_r}|,$$
setting $j_{r+1}=j_r+1$, we obtain condition (\ref{e2}) for $r$.

\bigskip

\noindent {\bf Case 2} {\it $I_{j_r}$ is clean and ${\bar{x}_t\,\,\bar{x}_{t+1}\choose x_t\,\,x_{t+1}}={0\,\,0\choose 0\,\,1}$.}

\bigskip 

\noindent Since $t\not\in \bar{I}$, we have $\bar{x}_{t-1}=1$, and hence, $x_{t-1}=0$.
$C(t,t+1)$ implies $x_{t+2}=1$.
Since $I_{j_r}$ is clean, $\bar{x}_{t+2}=0$.
$C(t+1,t+2)$ implies $x_{t+3}=1$.
Since $I_{j_r}$ is clean, $\bar{x}_{t+3}=0$.
This implies $i_{j_r+1}-i_{j_r}\geq 4$.
If $i_{j_r+1}-i_{j_r}\geq 5$, 
then setting $j_{r+1}=j_r+1$,
we obtain condition (\ref{e2}) for $r$ as in Case 1.
Hence, we may assume that $i_{j_r+1}-i_{j_r}=4$.

If $I_{j_r+1}$ is clean,
then $t+4\not\in \bar{I}$ implies $\bar{x}_{t+5}=1$, and hence, $x_{t+5}=0$.
Now, analogous arguments as in Case 1 imply
$x_{t+6}=x_{t+7}=x_{t+8}=0$
and
$\bar{x}_{t+6}=\bar{x}_{t+7}=\bar{x}_{t+8}=1$.
Hence, $i_{j_r+2}-i_{j_r}\geq 9$, and
$$\left|\left(C\cup \bar{C}\right)\cap V(I_{j_r}\cup I_{j_r+1})\right|
=|I_{j_r}\cup I_{j_r+1}|-2\geq \frac{7}{9}|I_{j_r}\cup I_{j_r+1}|,$$
that is, setting $j_{r+1}=j_r+2$, we obtain condition (\ref{e2}) for $r$.
Hence, we may assume that $I_{j_r+1}$ is dirty.

Let $j'_{r+1}$ be maximum such that $j_r<j'_{r+1}\leq k$ and
$I_j$ is dirty for $j\in \{ j_r+1,j_r+2,\ldots,j'_{r+1}-2\}$.
Clearly, $j'_{r+1}=k$ or $I_{j'_{r+1}-1}$ is clean.

If $j'_{r+1}=k$, then set $\ell=r+1$ and $j_{r+1}=k$.
Note that, if $I_{\ell-1}$ is dirty, then 
$\left|I_{j_r}\cup \cdots \cup I_{\ell-1}\right|\geq 6$ and 
$$
\left|
\left(C\cup \bar{C}\right)
\cap 
V(I_{j_r}\cup \cdots \cup I_{j_\ell-1})
\right|=\left|I_{j_r}\cup \cdots \cup I_{j_\ell-1}\right|-1
\geq \frac{7}{9}\left|I_{j_r}\cup \cdots \cup I_{j_\ell-1}\right|,$$
and, if $I_{j_\ell-1}$ is clean, then 
$\left|I_{j_r}\cup \cdots \cup I_{j_\ell-1}\right|\geq 8$ and 
$$
\left|
\left(C\cup \bar{C}\right)
\cap 
V(I_{j_r}\cup \cdots \cup I_{j_\ell-1})
\right|=\left|I_{j_r}\cup \cdots \cup I_{j_\ell-1}\right|-2
\geq \frac{7}{9}\left|I_{j_r}\cup \cdots \cup I_{j_\ell-1}\right|-\frac{5}{9},$$
that is, in both cases (\ref{e3}) holds.
Hence, we may assume that $j'_{r+1}<k$ and $I_{j'_{r+1}-1}$ is clean.

If $i_{j'_{r+1}-1}-i_{j_r+1}\geq 3$, 
then set $j_{r+1}=j'_{r+1}$.
Since 
$\left|I_{j_r}\cup \cdots \cup I_{j_{r+1}-1}\right|\geq 9$ and 
$$
\left|
\left(C\cup \bar{C}\right)
\cap 
V(I_{j_r}\cup \cdots \cup I_{j_{r+1}-1})
\right|=\left|I_{j_r}\cup \cdots \cup I_{j_{r+1}-1}\right|-2
\geq \frac{7}{9}\left|I_{j_r}\cup \cdots \cup I_{j_{r+1}-1}\right|,$$
(\ref{e2}) holds for $r$.
Hence, we may assume that $i_{j'_{r+1}-1}-i_{j_r+1}=2$, 
which implies that $I_{j_r+1}$ has exactly two elements, 
and $j'_{r+1}=j_r+3$, that is, $I_{j_r+2}$ is clean.

Let $s=i_{j_{r+1}'-1}$. 
Note that $s=t+6$.
If ${\bar{x}_s\,\,\bar{x}_{s+1}\choose x_s\,\,x_{s+1}}={0\,\,1\choose 0\,\,0}$,
then, $s+1\not\in \bar{I}$ implies $\bar{x}_{s+2}=1$, 
which implies that $|I_{i_r+2}|\geq 3$.
Again, setting $j_{r+1}=j'_{r+1}$
yields $\left|I_{j_r}\cup \cdots \cup I_{j_{r+1}-1}\right|\geq 9$,
and (\ref{e2}) for $r$ follows as above.
Hence, we may assume that 
${\bar{x}_s\,\,\bar{x}_{s+1}\choose x_s\,\,x_{s+1}}={0\,\,0\choose 0\,\,1}$, 
that is,
$$
{\bar{x}_t\,\,\ldots\,\,\bar{x}_{t+7}\choose x_t\,\,\ldots\,\,x_{t+7}}=
{
0\,\,0\,\,0\,\,0\,\,0\,\,1\,\,0\,\,0
\choose
0\,\,1\,\,1\,\,1\,\,0\,\,1\,\,0\,\,1
}.
$$
Now, $\bar{C}(t+4,t+6)$ does not hold, which is a contradiction,
and completes the second case.

\bigskip

\noindent For the remaining cases, we may assume that $I_{j_r}$ is dirty.
Let $j'_{r+1}$ be maximum such that $j_r<j'_{r+1}\leq k$ and
$I_j$ is dirty for $j\in \{ j_r,j_r+1,\ldots,j'_{r+1}-2\}$.
Clearly, $j'_{r+1}=k$ or $I_{j'_{r+1}-1}$ is clean.

If $j'_{r+1}=k$, then set $\ell=r+1$ and $j_{r+1}=k$.
Note that, if $I_{\ell-1}$ is dirty, then 
$$
\left|
\left(C\cup \bar{C}\right)
\cap 
V(I_{j_r}\cup \cdots \cup I_{j_\ell-1})
\right|=\left|I_{j_r}\cup \cdots \cup I_{j_\ell-1}\right|
\geq \frac{7}{9}\left|I_{j_r}\cup \cdots \cup I_{j_\ell-1}\right|,$$
and, if $I_{j_\ell-1}$ is clean, then 
$\left|I_{j_r}\cup \cdots \cup I_{j_\ell-1}\right|\geq 4$ and 
$$
\left|
\left(C\cup \bar{C}\right)
\cap 
V(I_{j_r}\cup \cdots \cup I_{j_\ell-1})
\right|=\left|I_{j_r}\cup \cdots \cup I_{j_\ell-1}\right|-1
\geq \frac{7}{9}\left|I_{j_r}\cup \cdots \cup I_{j_\ell-1}\right|-\frac{5}{9},$$
that is, in both cases (\ref{e3}) holds.
Hence, we may assume that $j'_{r+1}<k$ and $I_{j'_{r+1}-1}$ is clean.

The remaining two cases have some similarities with Cases 1 and 2.

Let $t=i_{j'_{r+1}-1}$.

\bigskip

\noindent {\bf Case 3} {\it $I_{j_r}$ is dirty and ${\bar{x}_t\,\,\bar{x}_{t+1}\choose x_t\,\,x_{t+1}}={0\,\,1\choose 0\,\,0}$.}

\bigskip 

\noindent Since $t\not\in \bar{I}$, we have $\bar{x}_{t+2}=1$.
Since $I_{j'_{r+1}-1}$ is clean, $x_{t+2}=0$.
This implies $i_{j'_{r+1}}-i_{j_r}\geq 5$.
Setting $j_{r+1}=j'_{r+1}$, condition (\ref{e2}) for $r$ follows as in Case 1.

\bigskip

\noindent {\bf Case 4} {\it $I_{j_r}$ is dirty and ${\bar{x}_t\,\,\bar{x}_{t+1}\choose x_t\,\,x_{t+1}}={0\,\,0\choose 0\,\,1}$.}

\bigskip 

\noindent If $i_{j'_{r+1}}-i_{j_r}\geq 5$, 
then setting $j_{r+1}=j'_{r+1}$ satisfies (\ref{e2}) for $r$ as in Case 3.
Hence, we may assume that $i_{j'_{r+1}}-i_{j_r}=4$,
which implies $j'_{r+1}=j_r+2$ and $|I_{j_r}|=|I_{j_r+1}|=2$. 

If $I_{j_r+2}$ is clean, 
then $t+2\not\in \bar{I}$ implies $\bar{x}_{t+3}=1$, and hence $x_{t+3}=0$.
Now similar arguments as in Case 1 imply
$x_{t+4}=x_{t+5}=x_{t+6}=0$ and
$\bar{x}_{t+4}=\bar{x}_{t+5}=\bar{x}_{t+6}=1$.
Therefore, $i_{j_r+3}-i_{j_r}\geq 9$, 
and setting $j_{r+1}=i_{j_r+3}$ satisfies (\ref{e2}) for $r$ as above.
Note that if $i_{j_r+3}-i_{j_r}=9$, then
$
{
\bar{x}_t\,\,\ldots\,\,\bar{x}_{t+8}
\choose 
x_t\,\,\ldots\,\,x_{t+8}
}
$
corresponds to the pattern used in the proof of Lemma \ref{lemma2}.
Hence, we may assume that $I_{j_r+2}$ is dirty.

Let $j''_{r+1}$ be maximum such that $j_r+2<j''_{r+1}\leq k$ and
$I_j$ is dirty for $j\in \{ j_r+2,j_r+3,\ldots,j''_{r+1}-2\}$.
Clearly, $j''_{r+1}=k$ or $I_{j''_{r+1}-1}$ is clean.
If $j''_{r+1}=k$, then setting $\ell=r+1$ and $j_{r+1}=k$,
and arguing similarly as in Case 2 yields (\ref{e3}).
Hence, we may assume $j''_{r+1}<k$ and $I_{j''_{r+1}-1}$ is clean.

If $i_{j''_{r+1}}-i_{j_r}\geq 9$, then setting $j_{r+1}=j''_{r+1}$ yields (\ref{e2}) for $r$ as above.
Hence, we may assume that $i_{j''_{r+1}}-i_{j_r}=8$,
which implies that $j''_{r+1}=j_r+4$ and $|I_{j_r+2}|=|I_{j_r+3}|=2$.
This implies 
$$
{\bar{x}_{t-2}\,\,\ldots\,\,\bar{x}_{t+6}\choose x_{t-2}\,\,\ldots\,\,x_{t+6}}\in
\left\{
{
0\,\,1\,\,0\,\,0\,\,0\,\,1\,\,0\,\,1\,\,0
\choose
0\,\,1\,\,0\,\,1\,\,0\,\,1\,\,0\,\,0\,\,0
},
{
0\,\,1\,\,0\,\,0\,\,0\,\,1\,\,0\,\,0\,\,0
\choose
0\,\,1\,\,0\,\,1\,\,0\,\,1\,\,0\,\,1\,\,0
}
\right\}.
$$
Now the first options leads to the contradiction $t+5\in \bar{I}$, and 
the second option leads to the contradiction that $\bar{C}(t+2,t+4)$ does not hold.

This completes the proof of Claim \ref{claim3}. $\Box$

\bigskip

\noindent If $k\leq 1$, then $n>15$ implies 
$$\left|\left(C\cup \bar{C}\right)\right|
\geq 
\left|\left(C\cup \bar{C}\right)\cap V([n-6]\setminus [7])
\right|\geq n-6-7-1=n-14>\frac{7}{9}n-12.$$
Hence, we may assume that $k\geq 2$.

Since $i_1$ is the smallest integer $i\geq 8$ with ${\bar{x}_i\choose x_i}={0\choose 0}$,
we have 
$$\left|
\left(C\cup \bar{C}\right)
\cap 
V([i_1-1])
\right|\geq i_1-1-7.$$
Since $i_k$ is the largest integer $i\leq n-6$ with ${\bar{x}_i\choose x_i}={0\choose 0}$,
we have 
$$\left|
\left(C\cup \bar{C}\right)
\cap 
V([n]\setminus [i_k-1])
\right|\geq n-6-i_k.$$
By Claim \ref{claim3}, we obtain
\begin{eqnarray*}
&& \left|
\left(C\cup \bar{C}\right)
\cap 
V([i_k-1]\setminus [i_1-1])
\right|\\
&=&\left|
\left(C\cup \bar{C}\right)
\cap 
V(I_1\cup \cdots \cup I_{k-1})
\right|\\
&=&
\left|
\left(C\cup \bar{C}\right)
\cap 
V(I_{j_1}\cup \cdots \cup I_{j_2-1})
\right|
+\sum_{r=2}^{\ell-2}
\left|
\left(C\cup \bar{C}\right)
\cap 
V(I_{j_r}\cup \cdots \cup I_{j_{r+1}-1})
\right|\\
&&+
\left|
\left(C\cup \bar{C}\right)
\cap 
V(I_{j_{\ell-1}}\cup \cdots \cup I_{j_\ell-1})
\right|\\
&\geq &
\left(\frac{7}{9}\left|I_{j_1}\cup \cdots \cup I_{j_2-1}\right|-\frac{5}{9}\right)
+\sum_{r=2}^{\ell-2}\left(\frac{7}{9}\left|I_{j_r}\cup \cdots \cup I_{j_{r+1}-1}\right|\right)
+\left(\frac{7}{9}\left|I_{j_{\ell-1}}\cup \cdots \cup I_{j_\ell-1}\right|-\frac{5}{9}\right)\\
&=& 
\frac{7}{9}\left|I_1\cup\cdots\cup I_{k-1}\right|
-\frac{10}{9}\\
&= & 
\frac{7}{9}(i_k-i_1)-\frac{10}{9}\\
&= & 
\frac{7}{9}\left((i_k-i_1)-\frac{10}{7}\right).
\end{eqnarray*}
Altogether, this implies
\begin{eqnarray*}
|C\cup \bar{C}| 
& \geq & 
(i_1-1-7)
+\left(\frac{7}{9}(i_k-i_1)-\frac{10}{9}\right)
+(n-6-i_k)\\
& \geq & 
\frac{7}{9}\left((i_1-1-7)+\left((i_k-i_1)-\frac{10}{7}\right)+(n-6-i_k)\right)\\
& \geq & \frac{7}{9}n-\frac{108}{9}\\
& = & \frac{7}{9}n-12,
\end{eqnarray*}
which completes the proof.
$\Box$

\bigskip

\noindent We proceed to our main result.

\begin{theorem}
${\rm ic}(C_n\bar{C}_n)=\frac{7}{9}n+\Theta(1)$ for $n\geq 3$.
\end{theorem}
{\it Proof:} This follows immediately from Lemma \ref{lemma2} and Lemma \ref{lemma4}. 
$\Box$


\begin{thebibliography}{}
\bibitem{a} D. Auger, Minimal identifying codes in trees and planar graphs with large girth, Eur. J. Comb. 31 (2010) 1372-1384.
\bibitem{bl} Y. Ben-Haim and S. Litsyn, Exact minimum density of codes identifying vertices in the square grid, SIAM J. Discrete Math. 19 (2005) 69-82.
\bibitem{bchl} N. Bertrand, I. Charon, O. Hudry, and A. Lobstein, 1-identifying codes on trees, Australas. J. Comb. 31 (2005) 21-35.
\bibitem{bcmms} M. Blidia, M. Chellali, F. Maffray, J. Moncel, and A. Semri, Locating-domination and identifying codes in trees, Australas. J. Comb. 39 (2007) 219-232.
\bibitem{chl} I. Charon, O. Hudry, and A. Lobstein, Minimizing the size of an identifying or locating-dominating code in a graph is NP-hard, Theor. Comput. Sci. 3 (2003) 2109-2120.
\bibitem{cghlmm} I. Charon, S. Gravier, O. Hudry, A. Lobstein, M. Mollard, and J. Moncel, A linear algorithm for minimum 1-identifying codes in oriented trees, Discrete Appl. Math. 154 (2006) 1246-1253.
\bibitem{cghlmpz} G. Cohen, S. Gravier, I. Honkala, A. Lobstein, M. Mollard, C. Payan, and G. Z\'{e}mor, Improved identifying codes for the grids, Electr. J. Comb. 6 (1) (1999) $\#$R19 (comment).
\bibitem{chlz} G. Cohen, I. Honkala, A. Lobstein, and G. Z\'{e}mor, New Bounds for Codes Identifying Vertices in Graphs, Electr. J. Comb. 6 (1) (1999) $\#$R19.
\bibitem{cmr} B. Courcelle, J.A. Makowsky, and U. Rotics, Linear Time Solvable Optimization Problems on Graphs of Bounded Clique-Width, Theory Comput. Systems 33 (2000) 125-150.
\bibitem{dgm} M. Daniel, S. Gravier, and J. Moncel, Identifying Codes in Some Subgraphs of the Square Lattice, Theor. Comput. Sci. 319 (2004) 411-421.
\bibitem{gw} W. Goddard and K. Wash, ID Codes in Cartesian Products of Cliques,
J. Combin. Math. Combin. Comput. 85 (2013) 97-106.
\bibitem{gms1} S. Gravier, J. Moncel, and A. Semri, Identifying codes of cycles, Eur. J. Comb. 27 (2006) 767-776.
\bibitem{gms2} S. Gravier, J. Moncel, and A. Semri, Identifying codes of Cartesian product of two cliques of the same size, Electr. J. Comb. 15 (2008) $\#$N4.
\bibitem{hik} R. Hammack, W. Imrich, and S. Klav\v{z}ar, Handbook of product graphs, 2nd ed. Discrete Mathematics and Its Applications, Boca Raton (2011).
\bibitem{hhsv} T.W. Haynes, M.A. Henning, P.J. Slater, and L.C. van der Merwe, The complementary product of two graphs, Bull. Inst. Comb. Appl. 51 (2007) 21-30.
\bibitem{hhks} T.W. Haynes, K.R.S. Holmes, D.R. Koessler, and L. Sewell, Locating-domination in complementary prisms of paths and cycles, Congr. Numerantium 199 (2009) 45-55.
\bibitem{jl} V. Junnila and T. Laihonen, Optimal lower bound for 2-identifying codes in the hexagonal grid, Electr. J. Comb. 19 (2012) $\#$P38.
\bibitem{kcl} M.G. Karpovsky, K. Chakrabarty, and L.B. Levitin, On a New Class of Codes for Identifying Vertices in Graphs, IEEE Transactions on Information Theory 44 (1998) 599-611.
\bibitem{lr} V. Lozin and D. Rautenbach, The relative clique-width of a graph, J. Combin. Theory Ser. B 97 (2007) 846-858.
\bibitem{ms} R. Martin and B. Stanton, Lower bounds for identifying codes in some infinite grids, Electr. J. Comb. 17 (2010) $\#$R122.
\end{thebibliography}
\end{document}